\documentclass[11pt]{article}

\newtheorem{thm}{Theorem}[section]

\newtheorem{lemma}[thm]{Lemma}

\newtheorem{cor}[thm]{Corollary}

\newtheorem{conj}[thm]{Conjecture}

\newtheorem{claim}[thm]{Claim}

\newcommand{\beq}[1]{\begin{equation}\label{#1}}
\newcommand{\enq}[0]{\end{equation}}

\newcommand{\qed}[0]{\begin{flushright} \rule{2mm}{3mm} \end{flushright}}

\usepackage{graphicx}

\begin{document}

\title{On homomorphisms from the Hamming cube to {\bf Z}}
\author{David Galvin\\
Microsoft Research\\
One Microsoft Way\\
Redmond, WA 98052} \maketitle

\renewcommand{\thefootnote}{\fnsymbol{footnote}}
 \footnotetext{Research
supported by a Graduate School Fellowship from Rutgers
University.}

\begin{abstract}
Write ${\cal F}$ for the set of homomorphisms from $\{0,1\}^d$ to
${\bf Z}$ which send $\underline{0}$ to $0$ (think of members of
${\cal F}$ as labellings of $\{0,1\}^d$ in which adjacent strings
get labels differing by exactly $1$), and ${\cal F}_i$ for those
which take on exactly $i$ values.
We give asymptotic formulae for
$|{\cal F}|$ and $|{\cal F}_i|$.

In particular, we show that the probability that a uniformly
chosen member ${\bf f}$ of ${\cal F}$ takes more than five values
tends to $0$ as $d \rightarrow \infty$.
This settles a conjecture
of J. Kahn.
Previously, Kahn had shown that there is a constant
$b$ such that ${\bf f}$ a.s. takes at most $b$ values.
This in
turn verified a conjecture of I. Benjamini {\em et al.}, that for
each $t > 0$, ${\bf f}$ a.s. takes at most $td$ values.

Determining $|{\cal F}|$ is equivalent both to counting the number
of rank functions on the Boolean lattice $2^{[d]}$ (functions $f
\colon 2^{[d]} \longrightarrow {\bf N}$ satisfying
$f(\emptyset)=0$ and $f(A) \leq f(A \cup x) \leq f(A)+1$ for all
$A \in 2^{[d]}$ and $x \in [d]$) and to counting the number of
proper $3$-colourings of the discrete cube (i.e., the number of
homomorphisms from $\{0,1\}^d$ to $K_3$, the complete graph on $3$
vertices).

Our proof uses the main lemma from Kahn's proof of constant range,
together with some combinatorial approximation techniques
introduced by A. Sapozhenko.
\end{abstract}

\section{Introduction}

\subsection{Background and statement of the result}

Write $Q_d$ for the $d$-dimensional Hamming cube (the graph whose
vertex set is $\{0,1\}^d$ and in which two vertices are joined by
an edge if they differ in exactly one coordinate).
Set
$$
{\cal F}=\{f \colon V(Q_d)\rightarrow {\bf Z} \colon
f(\underline{0})=0 \mbox{ and } u \sim v \Rightarrow
|f(u)-f(v)|=1\}.
$$
(That is, ${\cal F}$ is the set of graph homomorphisms from $Q_d$
to ${\bf Z}$, normalized to vanish at $\underline{0}$.)

In \cite{BenjaminiHaggstromMossel}, this set of functions is studied
from a probabilistic point of view, a motivating idea being that a
typical element of ${\cal F}$ should exhibit stronger concentration
behavior than an arbitrary element. Put uniform probability measure
on ${\cal F}$, and define the function $R$ on ${\cal F}$ by
$R(f)=\{f(v)\colon v\in V(Q_d)\}$ ($R$ is the {\bf range} of $f$).
In \cite{BenjaminiHaggstromMossel} the following conjecture is made
about the concentration of $|R|$:

\begin{conj}
For each $t>0$, ${\bf P}(|R|>td) \rightarrow 0$ as $d\rightarrow
\infty$.
\end{conj}

In \cite{Kahn}, something stronger is proved, and something
stronger still conjectured:

\begin{thm} \label{rangeconstant}
There is a constant $b$ such that ${\bf P}(|R|>b)=e^{-\Omega(d)}$.
\end{thm}

\begin{conj} \label{conj5isright}
${\bf P}(|R|>5)=e^{-\Omega(d)}$ and ${\bf P}(|R|=5)=\Omega(1)$.
\end{conj}

In this paper we prove Conjecture \ref{conj5isright} by
(asymptotically) counting the number of homomorphisms with various
ranges.
Specifically, if we set
$$
{\cal F}_i = \{f\in {\cal F}\colon |R(f)|=i\},
$$
we prove
\begin{thm} \label{ourresult}
\begin{eqnarray}
|{\cal F}| & = & (2e \pm e^{-\Omega(d)})2^{2^{d-1}} \nonumber \\
|{\cal F}_3| & = & (2 \pm e^{-\Omega(d)})2^{2^{d-1}} \nonumber \\
|{\cal F}_4| & = & (4\sqrt{e}-4 \pm e^{-\Omega(d)})2^{2^{d-1}} \nonumber \\
|{\cal F}_5| & = & (2e - 4\sqrt{e} +2 \pm e^{-\Omega(d)})2^{2^{d-1}}, \nonumber
\end{eqnarray}
\end{thm}
\noindent which gives Conjecture \ref{conj5isright}.
Setting
${\cal F}_{\leq 5} = \cup_{i \leq 5} {\cal F}_i$, we see that
Theorem \ref{ourresult} has the following weaker but more
elegantly formulated consequence:

\begin{cor} \label{nicecor}
$|{\cal F}| \sim |{\cal F}_{\leq 5}| \sim 2e2^{2^{d-1}}$.
\end{cor}

Corollary \ref{nicecor} makes sense: a little thought suggests
that a typical member of ${\cal F}$ should be constant on either
even or odd vertices of the cube, except for a small set of
``blemishes'' on which it takes values $2$ away from the
predominant value, and take just two values on
vertices of the other parity.

The problem under discussion is equivalent to the question of the
number of rank functions on the Boolean lattice $2^{[d]}$ (here
$[d]=\{1, \ldots, d\}$). A {\bf rank function} is an $f\colon
2^{[d]} \longrightarrow {\bf N}$ satisfying $f(\emptyset)=0$ and
$f(A) \leq f(A \cup x) \leq f(A)+1$ for all $A \in 2^{[d]}$ and $x
\in [d]$. An easy lower bound on the number of rank functions is
$2^{2^{d-1}}$ (consider those functions which take the value $k/2$
on each element of the $k$th level of the Boolean lattice for each
even $k$). Athanasiadis \cite{Athanasiadis} conjectured that the
total number of rank functions is $2^{2^{d-1}(1+o(1))}$. This
conjecture is proved in \cite{KahnLawrenz}, where it is further
conjectured that the number is in fact $O(2^{2^{d-1}})$. Theorem
\ref{ourresult} answers this conjecture in the affirmative; for, as
observed by Mossel (see \cite{Kahn}), there is a bijection from the
set of rank functions to ${\cal F}$: identifying a subset $A$ of
$[d]$ with a vertex of $Q_d$ in the natural way, the bijection is
given by $g \longrightarrow f$ where $f(A)=2g(A)-|A|$.

Theorem \ref{ourresult} also provides information about the number
of proper $3$-colourings of $Q_d$. A {\bf proper $3$-colouring} of a
graph $G$ with vertex set $V$ and edge set $E$ is a function
$\chi\colon V \longrightarrow \{0,1,2\}$ satisfying $(x,y) \in E
\Rightarrow \chi(x) \neq \chi (y)$. Theorem \ref{ourresult} implies
that the number of proper $3$-colourings of $Q_d$ is asymptotic to
$6e2^{2^{d-1}}$; for, as observed by Randall \cite{Randall2}, there
is a bijection from ${\cal F}$ to the set of proper $3$-colourings
of $Q_d$ with $\chi(\underline{0}) = 0$: the bijection is given by
$f \longrightarrow \chi$ where $\chi(v)=i$ iff $f(v) \equiv i$ (mod
$3$).

The main inspiration for the proof of Theorem \ref{ourresult} is the
work of A. Sapozhenko, who, in \cite{Sapozhenko2}, gave a relatively
simple derivation for the asymptotics of the number of independent
sets in $Q_d$ (earlier derived in a more involved way in
\cite{KorshunovSapozhenko}). Our Lemma \ref{mainhmsapprox} is a
modification of a lemma in \cite{Sapozhenko}, and our overall
approach is similar to \cite{Sapozhenko2}. The other key ingredient
in our proof is the main lemma from \cite{Kahn}, which was already
used by Kahn to give Theorem \ref{rangeconstant}.

In the rest of this section, we establish basic notation and
gather together the main external ingredients that will be used in
the proof of Theorem \ref{ourresult}, before giving an outline of the rest of the paper.

\subsection{Notation and conventions} \label{subsectionnotation}

For graph theory basics, see e.g. \cite{Bollobas},~\cite{Diestel}.
For basics of the combinatorics of the Hamming cube, see e.g.
\cite{Bollobas4}.

The Hamming cube $Q_d$ is a $d$-regular, bipartite graph.
Write
$V$ for the vertex set of the cube, ${\cal E}$ for the set of even
vertices (those whose $\ell_1$ distance from $\underline{0}$ is
even) and ${\cal O}$ for the set of odd vertices.
Set
$M=2^{d-1}=|{\cal E}|=|{\cal O}|$.

For $u, v \in V$ and $A,C \subseteq V$ we write $u \sim v$ if
there is an edge in $Q_d$ joining $u$ and $v$, $\nabla(A)$ for the
set of edges having exactly one end in $A$ and (when $A \cap
C=\emptyset$) $\nabla(A,C)$ for the set of edges having one end in
each of $A, C$.

Set $N(u)=\{w\in V\colon w \sim u\}$ ($N(u)$ is the {\bf
neighbourhood} of $u$), $N(A)=\cup_{w \in A} N(w)$, $N_C(u)=\{w\in
C\colon w \sim u\}$, $N_C(A)=\cup_{w \in A} N_C(w)$, and
$d_C(u)=|N_C(u)|$. Write $\rho(u,v)$ for the length of the
shortest $u$-$v$ path in $Q_d$, and set $\rho(u,A)=\min_{w \in
A}\{\rho(u,w)\}$ and $\rho(A,C)=\min_{w \in A, w' \in
C}\{\rho(w,w')\}$. Set $B(A)=\{v \in V\colon N(v) \subseteq A\}$.

We say that $A$ is {\bf $k$-linked} if for every $u,v\in A$ there
is a sequence $u=u_0, u_1, \ldots, u_l=v$ in $A$ with
$\rho(u_i,u_{i+1})\leq k$ for $i = 0, \ldots, l-1$.
Note that for
any $k$, $A$ is the disjoint union of its maximal $k$-linked
subsets --- we call these the {\bf $k$-components} of $A$.
Write
$C \prec A$ if $C$ is a $2$-component of $A$, and $c(A)$ for the
number of $2$-components of $A$.

We say that $A$ is {\bf small} if $|A|<\alpha^d$ for a certain
constant $\alpha < 2$ that will be discussed in Section
\ref{sectionreduction} (and {\bf large} otherwise), {\bf sparse}
if all the $2$-components of $A$ are singletons (and {\bf
non-sparse} otherwise), and {\bf nice} if $A$ is small, $2$-linked
and of size at least $2$.
Note that all sets $A$ that we will
consider will satisfy either $A \subseteq {\cal E}$ or $A
\subseteq {\cal O}$.

For integers $a < b$ we define $[a,b] = \{a, \ldots, b\}$.

We use ``$\ln$'' for the natural logarithm and ``$\log$'' always
means the base $2$ logarithm.
The implied constants in the $O$ and
$\Omega$ notation are absolute (independent of $d$).
We always
assume that $d$ is large enough to support our assertions.
No
attempt has been made to optimize constants.

\subsection{External ingredients}
\label{sectionext}

We list here the main results that we will be drawing on
in the rest of the paper.

\medskip

We begin with a lemma bounding the number of connected
subgraphs of a graph. The infinite $\Delta$-branching rooted tree contains precisely
${\Delta n \choose n}/((\Delta-1)n+1)$ rooted subtrees with $n$ vertices
(see e.g. Exercise 11 (p. 396) of \cite{Knuth})
and this implies that if $G$ is a graph with maximum degree $\Delta$ and vertex set $V(G)$
then the number of $n$-vertex subsets of $V(G)$ which contain a fixed vertex and induce a connected
subgraph is at most $(e\Delta)^{n}$. (This fact is rediscovered in \cite{Sapozhenko}.)
We will use the following easy corollary.

\begin{lemma} \label{Tree}
Let $\Sigma$ be a graph with vertex set $V(\Sigma)$ and
maximum degree $\Delta$. For each fixed $k$, the number of $k$-linked subsets of
$V(\Sigma)$ of size $n$ which contain a fixed
vertex is at most $2^{O(n\log \Delta)}$.
\end{lemma}

\noindent This follows from the fact that a $k$-linked subset of
$\Sigma$ is connected in a graph with all degrees
$O(\Delta^{k+1})$.

\medskip

The next lemma is a special case of a fundamental result due to
Lov\'asz \cite{Lovasz} and Stein \cite{Stein} (see also
\cite{Furedi}). For a bipartite graph $\Sigma$ with bipartition
$X\cup Y$, say $Y'\subseteq Y$ {\bf covers} $X$ if each $x\in X$
has a neighbour in $Y'$.

\begin{lemma}\label{Lcor}
If a bipartite graph $\Sigma$ with bipartition $X\cup Y$ satisfies
$d(x)\geq a$ for all $x\in X$ and $d(y)\leq b$ for all $y\in Y$,
then $X$ is covered by some $Y'\subseteq Y$ of size at most
$(|Y|/a)(1+\ln b)$.
\end{lemma}

\medskip

The next lemma is from \cite{Sapozhenko}
(see Lemma 2.1); the reader should have no
difficulty supplying a proof.

\begin{lemma} \label{Lconn}
If $\Sigma$ is a graph on vertex set $V(\Sigma)$ and $A,C \subseteq V(\Sigma)$ satisfy

\medskip
\noindent {\rm (i)}  $A$ is $k$-linked

\medskip
\noindent and

\medskip
\noindent {\rm (ii)}  $\rho(u,C)\leq l$ for each $u \in A$ and
$\rho(v,A)\leq l$ for each $v \in C$,

\medskip
\noindent then $C$ is $(k+2l)$-linked.
\end{lemma}

The main step from the proof of Theorem \ref{rangeconstant} in
\cite{Kahn} (obtained via entropy arguments) will also be used
here. For $f \in {\cal F}$, set $C(f)=\{v \in V\colon f|_{N(v)}
\mbox{ is constant}\}$.

\begin{lemma} \label{entropylemma}
For $u \sim v$ and ${\bf f}$ drawn uniformly from ${\cal F}$,
${\bf P}(|\{u,v\}\cap C({\bf f})|=1)=1-e^{-\Omega(d)}$.
\end{lemma}

Finally, we need to know something about isoperimetry in the cube.
A {\bf Hamming ball centered at $x_0$} in $Q_d$ is any set of
vertices $B$ satisfying
$$
\{u \in V\colon \rho(u,x_0) \leq k\} \subseteq B \subset \{u \in
V\colon \rho(u,x_0) \leq k+1\}
$$
for some $k<d$. An {\bf even} (resp. {\bf odd}) {\bf Hamming ball}
is a set of vertices of the form $B \cap {\cal E}$ (resp. $B \cap
{\cal O}$) for some Hamming ball $B$. We use the following result
of K\"orner and Wei \cite{KornerWei}.

\begin{lemma} \label{kornerandwei}
For every $C \subseteq {\cal E}$ (resp. ${\cal O}$) and $D \subseteq V$, there
exists an even (resp. odd) Hamming ball $C'$ and a set $D'$ such that
$|C'|=|C|$, $|D'|=|D|$ and $\rho(C',D') \geq \rho(C,D)$.
\end{lemma}

\subsection{Outline}

The rest of the paper is organized as follows.

In Section \ref{sectionreduction} we use Lemma \ref{entropylemma}
to reduce Theorem \ref{ourresult} to the problem of counting the
number of homomorphisms which are predominantly $0$ on ${\cal E}$.
The easy lower bounds on the number of homomorphisms which take on
four and five values are given in Section
\ref{sectionlowerbounds}. In Section \ref{sectionsums} we examine
a general type of sum over small subsets of ${\cal E}$ and
establish some of its properties. In Section \ref{sectionthesum}
we write down an explicit sum of the type examined in Section
\ref{sectionsums} for the number of homomorphisms which are
predominantly $0$ on ${\cal E}$. The rest of the paper is devoted
to estimating this sum. In Section \ref{sectionisoperimetry} we
establish lower bounds on the sizes of neighbourhoods of
single-parity sets in the cube. In Section
\ref{sectionmainapproximation} we arrive at the heart of the
matter, showing that the set of nice subsets of ${\cal E}$ can be
``well-approximated'' in a precise sense by members of a ``small''
collection; this allows us to swiftly complete the proof of
Theorem \ref{ourresult} in Section \ref{sectionprovingo(1)}. We
postpone a more detailed outline of the latter portion of the
argument until the beginning of Section
\ref{sectionmainapproximation}. Finally, in Section
\ref{sectionremarks}, we make some brief remarks on the proof and
possible extensions of the techniques used.

\section{Reduction to mostly constant} \label{sectionreduction}

We begin the proof of Theorem \ref{ourresult} by using Lemma
\ref{entropylemma} to reduce the problem to that of counting
homomorphisms which mainly take a single value on ${\cal E}$.

There is an inherent odd-even symmetry in the problem; we now reformulate
slightly to make use of this. Write
$$
{\cal A}= \{f\colon V \rightarrow {\bf Z} \colon  u \sim v
\Rightarrow |f(u)-f(v)|=1\}
$$
and write ${\cal B}$ for the quotient of ${\cal A}$ by the equivalence relation
$$
f \equiv g \iff f - g \mbox{ is constant on $V$}.
$$
For each $f \in {\cal A}$ write $[f]$ for the equivalence class of $f$ in ${\cal B}$.
Noting that $R$ is constant on equivalence classes, we may define
$$
{\cal B}_i = \{[f]\in {\cal B}\colon |R(f)|=i\}.
$$
Clearly $|{\cal B}_i| = |{\cal F}_i|$ for each $i$ (${\cal F}$ is a complete set
of representatives for ${\cal B}$).

For $f \in {\cal A}$, we say that $f$ is {\bf mostly constant on
${\cal E}$} if there is some $c$ such that $\{v \in {\cal E}\colon
f(v) \neq c\}$ is small (see Section \ref{subsectionnotation} for
the definition of small; the constant $\alpha$ in that definition
will be specified in the proof of Lemma \ref{mostlyconstant}), and
we define {\bf mostly constant on ${\cal O}$} analogously. These
definitions respect the equivalence relation, so we may define
$$
{\cal B}^{{\cal E}} = \{[f] \in {\cal B}\colon f \mbox{ is mostly
constant on ${\cal E}$}\}.
$$
Define ${\cal B}^{{\cal O}}$ analogously.
By symmetry, $|{\cal B}^{{\cal E}}|=
|{\cal B}^{{\cal O}}|$ (any automorphism of $Q_d$ that sends ${\cal E}$ to ${\cal O}$
induces a bijection between the two sets).

\begin{lemma} \label{negligibleintersection}
$$
|{\cal B}^{{\cal E}} \cap {\cal B}^{{\cal O}}| = e^{-\Omega(d)}|{\cal B}|.
$$
\end{lemma}

\noindent {\em Proof: }To specify an
$[f] \in {\cal B}^{{\cal E}} \cap {\cal B}^{{\cal O}}$ we first
specify the predominant values of the representative $f$ on
${\cal E}$ and ${\cal O}$. W.l.o.g. we may assume that the predominant value on ${\cal E}$
is $0$, and so the predominant value on ${\cal O}$ is one of $\pm 1$.
We then specify the small sets from ${\cal E}$
and ${\cal O}$ on which $f$ does not take the predominant values,
and finally the values of $f$ on these small sets. Noting that
once $f(v)$ has been specified for any $v \in V$ there are most $2d+1$
values that $f$ can take on any other vertex and that $2^M$ is a
trivial lower bound on $|{\cal B}|$, we get
\begin{eqnarray}
|{\cal B}^{{\cal E}} \cap {\cal B}^{{\cal O}}| & \leq &
2\sum_{i, j \leq \alpha^d} {M \choose i}{M \choose j}(2d+1)^{i+j} \nonumber \\
 & \leq & e^{-\Omega(d)}|{\cal B}|. \nonumber
\end{eqnarray}
\qed

\begin{lemma} \label{mostlyconstant}
$$
|{\cal B}| = (2 \pm e^{-\Omega(d)})|{\cal B}_{\cal E}|.
$$
\end{lemma}

\noindent {\em Proof: }For $f \in {\cal A}$, set $C(f)=\{v \in
V\colon f|_{N(v)} \mbox{ is constant}\}$ (extending the definition
given in Section \ref{sectionext}). We choose a uniform member
${\bf [f]}$ of ${\cal B}$ by choosing ${\bf f}$ uniformly from
${\cal F}$. For ${\bf [f]}$ and $u,v \in V$, let $Q_u$ be the
event $\{u \in C({\bf f})\}$, $Q_{\overline{u}}$ the complementary
event, $Q_{u\overline{v}} = Q_u \cap Q_{\overline{v}}$ and
$Q_{\overline{u}\overline{v}}= Q_{\overline{u}} \cap
Q_{\overline{v}}$. Write $K_u = K_u({\bf f})$ for the set of
vertices that can be reached from $u$ in $C({\bf f})$ via steps of
size exactly $2$, and let $Q_{uv}^*$ be the event $\{v \in K_u\}$.
(Note that if $f, g \in {\cal A}$ are equivalent then $C(f)=C(g)$,
so all these events are well defined.)

Let $u$ and $v$ be two vertices of the same parity. We claim that
$Q_{\overline{u}\overline{v}} \cup Q_{uv}^*$ occurs with
probability $1-e^{-\Omega(d)}$. For, let $ua_1a_2 \ldots a_{2k-1}v$ be a $u$-$v$ path of
length at most $d$ (the diameter of $Q_d$). Writing $a_0$ for $u$ and
$a_{2k}$ for $v$, we have
$$
Q_{\overline{u}\overline{v}} \cup
Q_{uv}^* \supseteq \cap_{i=0}^{2k-1}
(Q_{a_i\overline{a_{i+1}}} \cup Q_{\overline{a_i}a_{i+1}}).
$$
By Lemma \ref{entropylemma}, ${\bf P}(Q_{a_i\overline{a_{i+1}}}
\cup Q_{\overline{a_i}a_{i+1}}) = 1-e^{-\Omega(d)}$ for each $i$.
Hence ${\bf P}(Q_{\overline{u}\overline{v}} \cup Q_{uv}^*) \geq
1-de^{-\Omega(d)} = 1-e^{-\Omega(d)}$, as claimed.

We therefore have, for fixed $u \in V$ and any $v$ of the same parity as $u$,
${\bf P}(Q_{uv}^*|Q_u)>1-c^{-d}$, where $c>1$ is fixed. So, conditioning on $Q_u$, we have
$$
{\bf E}(|\{v\colon \rho(u,v) \mbox{ {\em even}},v \not \in K_u
\}|)\leq (2/c)^d,
$$
so that, by Markov's Inequality (with the constant $c'$ chosen so
that $2/c < c' < 2$),
\begin{equation} \label{conditioning}
{\bf P}(|K_u|<M-(c')^d | Q_u) \leq (2/cc')^d=e^{-\Omega(d)}.
\end{equation}

If $u \not \in C({\bf f})$, then $K_u({\bf f}) = \emptyset$, so that
${\bf P}(|K_u|<M-(c')^d | Q_{\overline{u}})=1$.
By symmetry, ${\bf P}(Q_{u \overline{v}})$ is the same for every adjacent $u$ and $v$,
and this together with Lemma \ref{entropylemma} gives
$1/2+e^{-\Omega(d)} > {\bf P}(Q_u), {\bf P}(Q_{\overline{u}}) > 1/2-e^{-\Omega(d)}$.
Combining these observations with (\ref{conditioning}), we get
$$
{\bf P}(|K_u|<M-(c')^d) \leq 1/2+e^{-\Omega(d)}.
$$
Noting that ${\bf f}$ is constant on the neighbourhood of $K_u$,
this says (taking $u$ to be any vertex in ${\cal O}$) that there
is a constant $\beta < 2$ such that
$$
{\bf P}({\bf f} \mbox{ is constant on a subset of ${\cal E}$ of size at least } M-\beta^d)
       >1/2-e^{-\Omega(d)}.
$$
Taking $\alpha = \beta$ in the definition of small, this says
$$
|{\cal B}^{{\cal E}}| \geq (1/2 - e^{-\Omega(d)})|{\cal B}|.
$$
The lemma now follows from Lemma \ref{negligibleintersection}. \qed

It is now convenient to choose as a complete set of representatives for
${\cal B}^{{\cal E}}$ the collection
$$
{\cal F}^{{\cal E}} = \{f \in {\cal A}\colon  {\cal E}\setminus
f^{-1}(0) \mbox{ is small}\}.
$$
Set
$$
{\cal F}^{{\cal E}}_i=\{f \in {\cal F}^{{\cal E}}\colon
|R(f)|=i\}.
$$
Noting that $|{\cal F}^{{\cal E}}_3| \geq 2^M$, we see that
Theorem \ref{ourresult} will now follow from
\begin{thm} \label{mainresult}
\begin{eqnarray}
|{\cal F}^{{\cal E}}| & \leq & (e+e^{-\Omega(d)})2^M \label{estimatingB} \\
|{\cal F}_4^{{\cal E}}| & \geq & (2\sqrt{e}-2 - e^{-\Omega(d)})2^M \label{lowerboundonf4} \\
|{\cal F}_5^{{\cal E}}| & \geq & (e-2\sqrt{e}+1-e^{-\Omega(d)})2^M \label{lowerboundonf5}.
\end{eqnarray}
\end{thm}
It is this that we proceed to prove.

\section{Lower bounds on $|{\cal F}^{{\cal E}}_4|$ and $|{\cal F}^{{\cal E}}_5|$}
\label{sectionlowerbounds}

The aim of this section is to prove (\ref{lowerboundonf4}) and (\ref{lowerboundonf5}).

With each sparse $A \subseteq {\cal E}$ of size at least $2$
we associate a subset ${\cal F}^{{\cal E}}_5(A) \subseteq {\cal F}^{{\cal E}}_5$
of size
$$
(2^{|A|}-2)2^{M-d|A|}= 2^MM^{-|A|}(1-2^{-|A|+1})
$$
consisting
of those $f \in {\cal F}^{{\cal E}}_5$ for which
$R(f)=[-2,2]$ and $f^{-1}(\{\pm 2\})=A$
(on $A$, choose values for $f$ from $\{\pm 2\}$,
choosing at least one $2$ and at least one $-2$; on ${\cal E}\setminus A$ give $f$ value 0;
and on ${\cal O}\setminus N(A)$ choose values from
$\{\pm 1\}$, all choices made independently).
Then ${\cal F}^{{\cal E}}_5(A) \cap {\cal F}^{{\cal E}}_5(B) = \emptyset$ whenever $A \neq B$.
Noting that there are at least ${M \choose k}-Md^2{M-2 \choose k-2}$
sparse subsets of ${\cal E}$ of size $k$, and that for $k \leq d$, this number is
$(1-e^{-\Omega(d)}){M \choose k}$, we can lower bound $|{\cal F}^{{\cal E}}_5|$ by
\begin{eqnarray}
\left|{\cal F}^{{\cal E}}_5 \right| & \geq &
        2^M \sum_{k \geq 2} |\{A \subseteq {\cal E}\colon  A \mbox{ sparse, }
                 |A|=k\}|M^{-k}(1-2^{-k+1}) \nonumber \\
 & \geq & 2^M (1-e^{-\Omega(d)}) \sum_{k=2}^d {M \choose k}M^{-k}(1-2^{-k+1}) \nonumber \\
 & \geq & 2^M (1-e^{-\Omega(d)}) \sum_{k=2}^d (1-e^{-\Omega(d)})(1/k!)(1-2^{-k+1}) \nonumber \\
 & \geq & 2^M (1-e^{-\Omega(d)}) (\sum_{k=2}^d 1/k! - 2\sum_{k=2}^d 2^{-k}/k!) \nonumber \\
 & \geq & 2^M (1-e^{-\Omega(d)})((e-2)-2(\sqrt{e}-3/2)) \nonumber \\
 & \geq & 2^M (e-2\sqrt{e}+1-e^{-\Omega(d)}), \nonumber
\end{eqnarray}
so we have (\ref{lowerboundonf5}).

We do something
similar for (\ref{lowerboundonf4}). With each nonempty, sparse $A \subseteq {\cal E}$
we associate a subset ${\cal F}^{{\cal E}}_4(A) \subseteq {\cal F}^{{\cal E}}_4$
of size
$$
2^{1+M-d|A|} = 2^MM^{-|A|}2^{-|A|+1}
$$
consisting
of those $f \in {\cal F}^{{\cal E}}_4$ for which either
$R(f)=[-2,1]$ or $R(f)=[-1,2]$ and $f^{-1}(\{\pm 2\})=A$
(choose a value from $\pm 2$ for $f$ to take on $A$; on ${\cal E}\setminus A$ give $f$ value 0;
and choose values from
$\pm 1$ on ${\cal O}\setminus N(A)$, all choices made independently).
So we have
\begin{eqnarray}
\left|{\cal F}^{{\cal E}}_4\right| & \geq &
        2^M \sum_{k \geq 1} |\{A \subseteq {\cal E}\colon  A \mbox{ sparse, }
                 |A|=k\}|M^{-k}2^{-k+1} \nonumber \\
 & \geq & 2^M (2\sqrt{e}-2-e^{-\Omega(d)}). \nonumber
\end{eqnarray}

\section{Sums over small subsets of ${\cal E}$} \label{sectionsums}

In this section, we examine a certain kind of sum that will arise when we try
to write down an explicit expression for $|{\cal F}^{{\cal E}}|$. Specifically, we prove

\begin{lemma} \label{generalsums}
Suppose that $g\colon 2^{\cal E}\rightarrow{\bf R}^+$ satisfies
\beq{sumcond1} g(A) = \prod \{ g(A_i) \colon  A_i \prec A \}, \enq
\beq{sumcond2} g(\{y\}) = c2^{-d} ~\mbox{for all $y \in {\cal E}$
for some constant $c>0$} \enq and \beq{sumcond3} \sum_{\mbox{$A$
nice}} g(A) = e^{-\Omega(d)}. \enq Then for all $D \subseteq {\cal
E}$
$$
\left|\sum_{\mbox{$A \subseteq D$,~$A$ small}} g(A) -
(1+c2^{-d})^{|D|}\right| =  e^{-\Omega(d)}.
$$
\end{lemma}

\noindent {\em Remark: }Because $\emptyset \prec \emptyset$, any
$g$ satisfying (\ref{sumcond1}) must also satisfy
$g(\emptyset)=1$.

\medskip

\noindent {\em Proof of Lemma \ref{generalsums}: }All summations
below are restricted to subsets of $D$. We begin by observing that
$(1+c2^{-d})^{|D|}=\sum_A c^{|A|}2^{-d|A|}$ and that if $A$ is
sparse then $g(A)=c^{|A|}2^{-d|A|}$, so that
\begin{equation} \label{whee}
\left|\sum_{\mbox{$A$ small}} g(A)-(1+c2^{-d})^{|D|}\right| \leq
   \sum ~\hspace{-1.8mm}' g(A) + \sum ~\hspace{-1.8mm}'' c^{|A|}2^{-d|A|} +
   \sum ~\hspace{-1.8mm}''' c^{|A|}2^{-d|A|},
\end{equation}
where $\sum'$ is over $A$ small and non-sparse, $\sum''$ is over $A$ large and
$\sum'''$ is over $A$ non-sparse.

We bound each of the terms on the right-hand side of (\ref{whee}). For the first we have
\begin{eqnarray}
\sum ~\hspace{-1.8mm}' g(A)
            & \leq & \sum \left\{ g(A')g(A' \setminus A)\colon
                                   A' \mbox{ nice, }A \mbox{ small, }A' \prec A \right\} \nonumber \\
            & \leq & \sum_{\mbox{$A'$ nice}} g(A') \sum_{\mbox{$A$ small}} g(A) \nonumber \\
            & = & e^{-\Omega(d)}\sum_{\mbox{$A$ small}} g(A) \label{whee1}.
\end{eqnarray}

\noindent For the second we have
\begin{eqnarray}
\sum ~\hspace{-1.8mm}'' c^{|A|}2^{-d|A|} & \leq &
\sum_{|A| \geq d} c^{|A|}2^{-d|A|} \nonumber \\
& \leq & \sum_{i = d}^{|D|} {|D| \choose i} (c2^{-d})^i \nonumber \\
& \leq & \sum_{i \geq d} c^i/i! \nonumber \\
& = & e^{-\Omega(d)}. \label{whee2}
\end{eqnarray}

\noindent Finally, for the third we have
\begin{eqnarray}
\sum ~\hspace{-1.8mm}''' c^{|A|}2^{-d|A|} & \leq
        & \sum_{x,x' \in D,~\rho(x,x')=2} c^22^{-2d} \sum_{A} c^{|A|}2^{-d|A|} \nonumber \\
& \leq & |D|c^2d^22^{-2d}(1+c2^{-d})^{|D|} \nonumber \\
& = & e^{-\Omega(d)}. \label{whee3}
\end{eqnarray}

Combining (\ref{whee1}), (\ref{whee2}) and
(\ref{whee3}) we get
\begin{eqnarray}
\left|\sum_{\mbox{$A$ small}} g(A)-(1+c2^{-d})^{|D|}\right| & = &
e^{-\Omega(d)} \left(\sum_{\mbox{$A$ small}} g(A)+1 \right) \label{intone} \\
& = & e^{-\Omega(d)}. \label{inttwo}
\end{eqnarray}
(We get (\ref{inttwo}) from (\ref{intone}) because the latter
implies that $\sum_{\mbox{$A$ small}} g(A)$ is bounded.) \qed

The most important $g$ that we will be considering is
$$
g(A)= 2^{-|N(A)|+|B(A)|}
$$
(recall that $B(A)=\{v\in N(A)\colon N(v)\subseteq A\}$). It's
easy to see that this satisfies (\ref{sumcond1}) and
(\ref{sumcond2}) (with $c=1$). It is far from obvious that it
satisfies (\ref{sumcond3}); Sections
\ref{sectionmainapproximation} and \ref{sectionprovingo(1)} are
devoted to the proof of this fact, which we state now for use in
Section \ref{sectionthesum}.

\begin{thm} \label{heartofmatter}
$$
\sum_{\mbox{$A \subseteq {\cal E}$ nice}} 2^{-|N(A)|+|B(A)|} =
e^{-\Omega(d)}.
$$
\end{thm}

\section{Proof of (\ref{estimatingB})} \label{sectionthesum}

In this section, we write an explicit sum of the type introduced in Section \ref{sectionsums}
for $|{\cal F}^{{\cal E}}|$ and use Lemma \ref{generalsums} to estimate it, modulo Theorem \ref{heartofmatter}.
This will give (\ref{estimatingB}).

For each small $A \subseteq {\cal E}$, set
$$
{\cal F}^{{\cal E}}(A) = \{f \in {\cal F}^{{\cal E}}\colon
f^{-1}(0)={\cal E}\setminus A\}.
$$
We may specify an $f \in {\cal F}^{{\cal E}}(A)$ by the following
procedure. First, noting that $f$ must be either always positive
or always negative on a $2$-component of $A$, we specify a sign
($\pm$) for each such $2$-component. Next, we specify a nested
sequence
$$
A = C_2 \supseteq C_4 \supseteq \ldots \supseteq C_{2[d/2]}.
$$
For each $i=1, \ldots, [d/2]$, $C_{2i}=\{u \in {\cal E}\colon
|f(u)| \geq 2i\}$. Because the diameter of $Q_d$ is $d$, we have
$|f(u)| \leq 2[d/2]$ for all $u \in {\cal E}$, so this second step
completes the specification of $f$ on ${\cal E}$. Note that not
every sequence of $C_{2i}$'s gives rise to a legitimate $f \in
{\cal F}^{{\cal E}}$.

\begin{figure}[htb]
  \begin{center}
    \includegraphics[width=12cm]{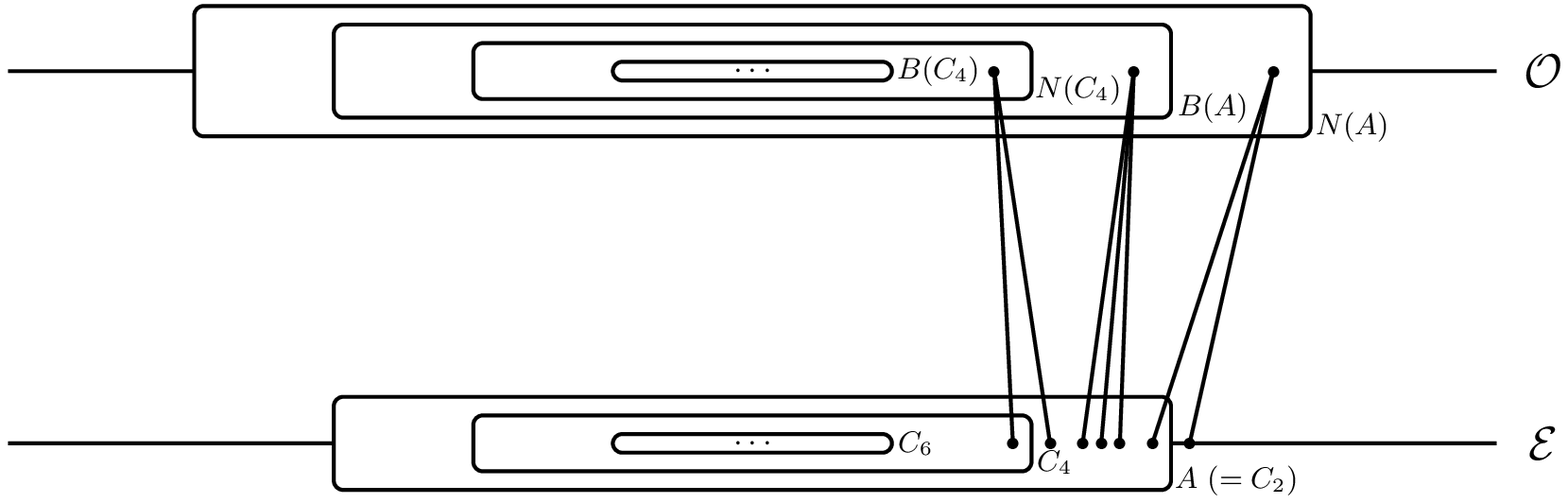}
    \caption{A vertex in $N(A)\setminus B(A)$ has neighbours in
    both ${\cal E}\setminus A$ and $A$, and a vertex in $N(C_4)\setminus
    B(C_4)$ has neighbours in both $A\setminus C_4$ and $C_4$, but
    a vertex in $B(A)\setminus N(C_4)$ only has neighbours in
    $A\setminus C_4$.}
      \label{figure1}
  \end{center}
\end{figure}

To specify $f$ on ${\cal O}$, we first specify a value from $\pm
1$ on each vertex of ${\cal O} \setminus N(A)$, and then, for each
$i=1, \ldots, [d/2]$, specify a value from $2i \pm 1$ for $|f(u)|$
for each $u \in B(C_{2i}) \setminus N(C_{2i+2})$ (note that the
sign of $f(u)$ for such $u$ has been determined by the
specification of signs on $A$). To see that this completes the
specification of $f$ on ${\cal O}$, note that we have a choice for
the value of $|f|$ at $u \in N(A)$ iff $f$ is constant on $N(u)$
iff $u \in B(C_{2i}) \setminus N(C_{2i+2})$ for some $1 \leq i
\leq [d/2]$ (setting $C_{2[d/2]+2}=\emptyset$), and that in this
case we can choose from two possible values, $2i \pm 1$ (see
Figure \ref{figure1}).

So, noting that $N(C_{2i+2}) \subseteq B(C_{2i})$ for each $i=1, \ldots, [d/2]$, we have
$$
|{\cal F}^{{\cal E}}(A)| = 2^{c(A)+ M-|N(A)|+|B(A)|} \sum \prod_{i=2}^{[d/2]} 2^{-|N(C_{2i})|+|B(C_{2i})|}
$$
where the sum --- here and in the next line --- is
over all legitimate choices of $C_2 \supseteq \ldots \supseteq C_{2[d/2]}$. Setting
$$
h(A)=2^{c(A)-|N(A)|+|B(A)|} \sum \prod_{i=2}^{[d/2]} 2^{-|N(C_{2i})|+|B(C_{2i})|}
$$
we get
$$
|{\cal F}^{{\cal E}}|=2^M \sum_{\mbox{$A \subseteq {\cal E}$ small}} h(A).
$$

We claim that $h$ satisfies all the conditions of Lemma
\ref{generalsums}. For $A=\{y\}$ we have $B(A)=\emptyset$, and so
$h(A)=2^{1-d}$; this gives (\ref{sumcond2}) (with $c=2$). To see
that $h$ satisfies (\ref{sumcond3}), note that for each $A
\subseteq {\cal E}$ small, each $C_{2i}$ is a small subset of $A$,
and so we can crudely upper bound $h(A)$ by
\begin{eqnarray}
h(A) & \leq & 2^{c(A)-|N(A)|+|B(A)|}
    \left( \sum_{\mbox{$C \subseteq A$ small}} 2^{-|N(C)|+|B(C)|} \right)^{[d/2]}\nonumber \\
 & \leq & 2^{c(A)-|N(A)|+|B(A)|}
    \left( (1+2^{-d})^{\alpha^d} + e^{-\Omega(d)} \right)^{[d/2]} \label{applyingsumsthm} \\
 & \leq & \left(1+o(1)\right)2^{c(A)-|N(A)|+|B(A)|}. \nonumber
\end{eqnarray}
The inequality in (\ref{applyingsumsthm}) is obtained by applying
Lemma \ref{generalsums} and Theorem \ref{heartofmatter}, and
(\ref{sumcond3}) for $h$ now follows directly from Theorem
\ref{heartofmatter}. Finally, to establish (\ref{sumcond1}) for
$h$, note that $C_2 \supseteq C_4 \supseteq \ldots \supseteq
C_{2[d/2]}$ is a legitimate sequence of $C$'s for $A$ iff $C_2\cap
A_i \supseteq C_4\cap A_i \supseteq \ldots \supseteq
C_{2[d/2]}\cap A_i$ is a legitimate sequence for $A_i$ for each
$2$-component $A_i$ of $A$, from which the claimed factorization
of $h(A)$ follows.

We can now easily establish (\ref{estimatingB}), thus completing
the proofs of Theorems \ref{mainresult} and \ref{ourresult}.
Applying Lemma \ref{generalsums}, we have
\begin{eqnarray}
\left||{\cal F}^{{\cal E}}| -e2^M\right| & \leq & 2^M\left(\left|
\sum ~\hspace{-1.8mm}' h(A)
-(1-2^{-d+1})^{|{\cal E}|}\right|
+\left|(1-2^{-d+1})^{|{\cal E}|}-e\right|\right) \nonumber \\
& = & e^{-\Omega(d)}2^{M}, \nonumber
\end{eqnarray}
\noindent where $\sum'$ is over $A \subseteq {\cal E}$ small.

\section{Isoperimetry in the cube} \label{sectionisoperimetry}

The aim of this section is to put some lower bounds on the
neighbourhood size of a small set in $Q_d$. We begin with

\begin{lemma} \label{boundsondelta}
For all $A \subseteq {\cal E}$ or $A \subseteq {\cal O}$ small, $|A|\leq
(1-\Omega(1))|N(A)|$.
\end{lemma}

\noindent {\em Proof: }By symmetry, we need only prove this when $A \subseteq {\cal E}$.
Let small $A \subseteq {\cal E}$ be given. Applying
Lemma \ref{kornerandwei} with $C=A$ and $D=V \setminus (A
\cup N(A))$, we find that there exists an even Hamming ball $A'$
with $|A'|=|A|$ and $|N(A)| \geq |N(A')|$. So we may assume that
$A$ is a small even Hamming ball.

We consider only the case where $A$ is centered at an even vertex,
w.l.o.g. $\underline{0}$, the other case being similar. In this
case,
$$
\{v \in {\cal E}\colon \rho(v,\underline{0}) \leq k\} \subseteq A
\subset \{v \in {\cal E}\colon \rho(v,\underline{0}) \leq k+2\}
$$
for some even $k \leq d/2-\Omega(d)$ (the bound on $k$ coming from
the fact that $A$ is small). For each $0 \leq i \leq (k+2)/2$, set
$B_i=A \cap \{v\colon \rho(v,\underline{0})=2i\}$, and
$N^+(B_i)=N(B_i)\cap\{u\colon \rho(u,\underline{0})=2i+1\}$. It's
clear that $N(A)=\cup_{0 \leq i \leq (k+2)/2} N^+(B_i)$ and that
for $i=0,\ldots,(k+2)/2$
\begin{eqnarray}
\frac{|B_i|}{|N^+(B_i)|} & \leq & \frac{2i+1}{d-2i} \label{LYM} \\
& = & 1-\Omega(1), \label{usingk}
\end{eqnarray}
\noindent from which the lemma follows. The inequality in
(\ref{usingk}) comes from the bound on $k$. The inequality in
(\ref{LYM}) is actually an equality except when $i=(k+2)/2$, in
which case it follows from the observation that each vertex in
$B_{k+2}$ has exactly $d-(k+2)$ neighbours in $N^+(B_{k+2})$, and
each vertex in $N^+(B_{k+2})$ has at most $(k+2)+1$ neighbours in
$B_{k+2}$. \qed

Lemma \ref{boundsondelta} is true for all small $A$, but can be strengthened
considerably when we impose stronger bounds on $|A|$. In this direction, we only need
the simple

\begin{lemma} \label{boundsondeltaAsmall}
If $|A|<d^{O(1)}$, then $|A| \leq O(1/d)|N(A)|$, and if $|A| \leq d/2$,
then $|N(A)| \geq d|A|-2|A|(|A|-1)$.
\end{lemma}

\noindent {\em Remark: }Note that the second statement is true for
all $A$, but vacuously so for $|A|>d/2$.

\medskip

\noindent {\em Proof of Lemma \ref{boundsondeltaAsmall}: }If $|A|<d^{O(1)}$,
then we have $k =O(1)$ in the notation of
Lemma \ref{boundsondelta}, and repeating the argument of that lemma
we get $|A| \leq O(1/d)|N(A)|$.

For the second part, note that each $u \in A$ has $d$ neighbours,
of which at least $d-2(|A|-1)$ must be unique to it, since a pair
of vertices in the cube can have at most two common neighbours.
\qed

From here on, the only properties of the cube that we will use are
the isoperimetric bounds of Lemmas \ref{boundsondelta} and
\ref{boundsondeltaAsmall}.

\section{The main approximation} \label{sectionmainapproximation}

We now begin the proof of Theorem \ref{heartofmatter}. The
approach will be to partition the set of $A$'s over which we are
summing according to the sizes of $A$, $N(A)$, $B(A)$ and
$N(B(A))$ (note that the summand in Theorem \ref{heartofmatter} is
constant on each partition class). The bulk of the work will be in
bounding the sizes of the partition classes.

Given $A \subseteq {\cal E}$, set $G=G(A)=N(A)$, $B=B(A)$ and
$H=H(A)=N(B)$. In what follows, $G$, $B$ and $H$ are always
understood to be $G(A)$, $B(A)$ and $H(A)$ for whatever $A$ is
under discussion. Note that $B \subseteq G$ and $H \subseteq A$.

Given $a$, $g$, $b$ and $h$, set
$$
{\cal H}(a, g, b, h)=\left\{A \subseteq {\cal E}~\mbox{ $2$-linked
$\colon |A|=a, |G|=g, |B|=b$ and $|H|=h$} \right\}.
$$
The aim of this section is to prove
\begin{lemma} \label{mainapprox}
For each $a$, $g$, $b$ and $h$ with $a \leq \alpha^d$,
$$
|{\cal H}(a,g,b,h)| < M 2^{g-b-\Omega(g/\log d)},
$$
\end{lemma}
from which we will
easily derive Theorem \ref{heartofmatter} in Section \ref{sectionprovingo(1)}.

\medskip

From now until the beginning of Section \ref{sectionprovingo(1)},
$a, g, b$ and $h$ are fixed, and we write ${\cal H}$ for ${\cal
H}(a,g,b,h)$. The proof of Lemma \ref{mainapprox} involves the
idea of ``approximation''. We begin with an informal outline. To
bound $|{\cal H}|$, we produce a small set ${\cal U}$ with the
properties that each $A \in {\cal H}$ is ``approximated'' (in an
appropriate sense) by some $U \in {\cal U}$, and for each $U \in
{\cal U}$, the number of $A \in {\cal H}$ that could possibly be
``approximated'' by $U$ is small. (Each $U \in {\cal U}$ will
consist of four parts; one each approximating $G$, $A$, $H$ and
$B$.) The product of the bound on $|{\cal U}|$ and the bound on
the number of $A \in {\cal H}$ that may be approximated by any $U$
is then a bound on $|{\cal H}|$. Another way of saying this is
that we produce a set ${\cal U}$ and a map $app\colon {\cal H}
\rightarrow {\cal U}$; we then bound $|{\cal H}|$ by
$$
|{\cal H}| \leq |{\cal U}|\max_{U \in {\cal U}}|app^{-1}(U)|.
$$
The set ${\cal U}$ is itself produced by an approximation process
--- we first produce a small set ${\cal V}$ with the property that
each $A \in {\cal H}$ is ``weakly approximated'' (in an
appropriate sense) by some $V \in {\cal V}$, and then show that
for each $V$ there is a small set ${\cal W}(V)$ with the property
that for each $A \in {\cal H}$ that is ``weakly approximated'' by
$V$, there is a $W \in {\cal W}(V)$ which approximates $A$; we
then take ${\cal U}=\cup_{V \in {\cal V}} {\cal W}(V)$. (Each $V
\in {\cal V}$ will consist of two parts; one each approximating
$G$ and $H$.)

\begin{figure}[htb]
  \begin{center}
    \includegraphics[width=12cm]{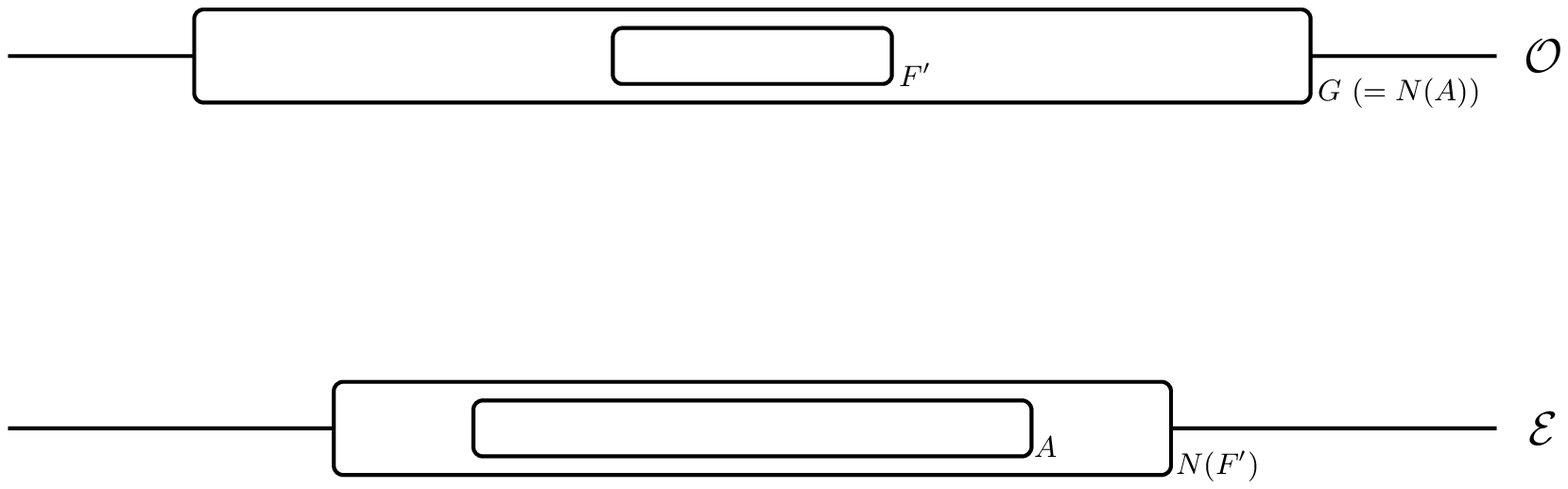}
    \caption{$F'$ satisfies both the conditions of (\ref{cover1}).}
      \label{figure2}
  \end{center}
\end{figure}

We now begin the formal discussion of Lemma \ref{mainapprox} by
introducing the two notions of approximation that we will use,
beginning with the weaker notion. A {\bf covering approximation}
for $A \subseteq {\cal E}$ is a pair $(F',P') \in 2^{\cal O}
\times 2^{\cal E}$ satisfying
\begin{equation} \label{cover1}
F' \subseteq G,~N(F') \supseteq A
\end{equation}
and
$$
P' \subseteq H,~N(P') \supseteq B
$$
\noindent (see Figure \ref{figure2}). An {\bf approximating
quadruple} for $A \subseteq {\cal E}$ is a quadruple $(F,S,P,Q)
\in 2^{\cal O} \times 2^{\cal E} \times 2^{\cal E} \times 2^{\cal
O}$ satisfying \beq{quad1} F \subseteq G,~S \supseteq A, \enq
\beq{quad2} d_F(u)>d-\sqrt{d}~~\mbox{for all $u \in S$} \enq
\beq{quad3} d_{{\cal E} \setminus S}(v) > d-\sqrt{d}~~\mbox{for
all $v \in {\cal O} \setminus F$} \enq \beq{quad4} P \subseteq
H,~Q \supseteq B, \enq \beq{quad5} d_P(u)>d-\sqrt{d}~~\mbox{for
all $u \in Q$} \enq and \beq{quad6} d_{{\cal O} \setminus Q}(v) >
d-\sqrt{d}~~\mbox{for all $v \in {\cal E} \setminus P$} \enq

\begin{figure}[htb]
  \begin{center}
    \includegraphics[width=12cm]{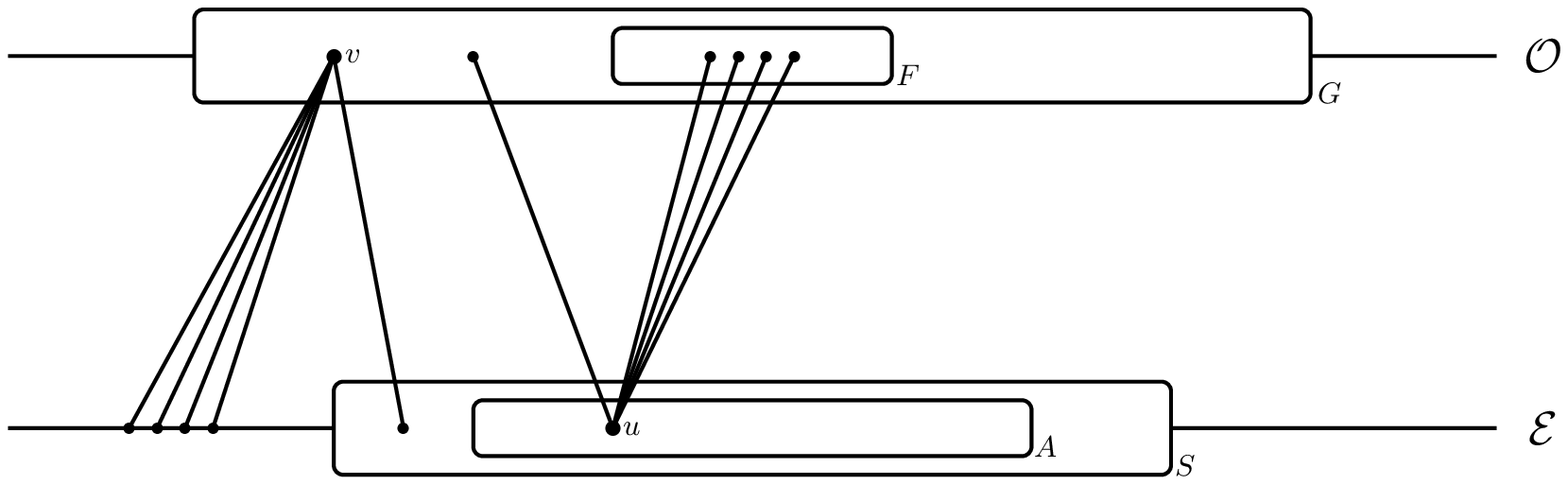}
    \caption{The pair $(F,S)$ satisfies (\ref{quad1}). To satisfy
    (\ref{quad2}) and (\ref{quad3}), each vertex $u \in S$ should have
    most (all but $\sqrt{d}$) of its neighbours in $F$, and each vertex $v \in {\cal O}\setminus F$
    should have most of its neighbours in ${\cal E}\setminus S$.}
      \label{figure3}
  \end{center}
\end{figure}

\noindent (see Figure \ref{figure3}). Note that if $x$ is in $A$
then all of its neighbours are in $G$, and if $y$ is in ${\cal O}
\setminus G$ then all of its neighbours are in ${\cal E} \setminus
A$. If we think of $S$ as ``approximate $A$'' and $F$ as
``approximate $G$'', (\ref{quad2}) says that if $x \in {\cal E}$
is in ``approximate $A$'' then almost all of its neighbours are in
``approximate $G$'', while (\ref{quad3}) says that if $y\in {\cal
O}$ is not in ``approximate $G$'' then almost all of its
neighbours are not in ``approximate $A$'', and there are similar
interpretations for (\ref{quad5}) and (\ref{quad6}).

\medskip

There are two parts to the proof of Lemma \ref{mainapprox}; the
``approximation'' step (Lemma \ref{mainhmsapprox}) and the
``reconstruction'' step (Lemma \ref{countingas0}). We now state
these two lemmas (from which Lemma \ref{mainapprox} follows
immediately).

\begin{lemma} \label{mainhmsapprox}
There is a family
$$
{\cal U}={\cal U}(a, g, b, h) \subseteq
2^{\cal O} \times 2^{\cal E} \times 2^{\cal E} \times 2^{\cal O}
$$
with
$$
|{\cal U}| \leq M2^{O(g\log d/\sqrt{d})}
$$
such that every $A \in {\cal H}$ has an approximating quadruple in
${\cal U}$.
\end{lemma}

\begin{lemma} \label{countingas0}
For each $(F,S,P,Q) \in 2^{{\cal O}} \times 2^{{\cal E}} \times
2^{{\cal E}} \times 2^{{\cal O}}$ satisfying (\ref{quad2}),
(\ref{quad3}), (\ref{quad5}) and (\ref{quad6}), there are at most
$2^{g - b - \Omega(g/\log d)}$ $A$'s in ${\cal H}$ satisfying
(\ref{quad1}) and (\ref{quad4}).
\end{lemma}

Lemma \ref{mainhmsapprox} follows directly from the next two
lemmas.

\begin{lemma} \label{firsthmsapprox}
There is a family
$$
{\cal V}={\cal V}(a,g,b,h) \subseteq 2^{\cal O} \times 2^{\cal E}
$$
with
$$
|{\cal V}| \leq M2^{O(g \log^2 d/d)}
$$
such that each $A \in {\cal H}$ has a covering approximation in
${\cal V}$.
\end{lemma}

\begin{lemma} \label{secondhmsapprox}
For each $(F',P') \in 2^{\cal O} \times 2^{\cal E}$ there is a
family
$$
{\cal W}={\cal W}(F',P',a,g,b,h) \subseteq 2^{\cal O}
\times 2^{\cal E} \times 2^{\cal E} \times 2^{\cal O}
$$
with
$$
|{\cal W}| \leq 2^{O(g \log d / \sqrt{d})}
$$
such that any $A \in {\cal H}$ for which $(F',P')$ is a covering
approximation has an approximating quadruple in ${\cal W}$.
\end{lemma}

We prove Lemmas \ref{firsthmsapprox} and \ref{secondhmsapprox} in
Section \ref{subsectionapproxlemmas}. We then prove Lemma
\ref{countingas0} in Section \ref{subsectionprovingmainapprox}.
The main point in the proof of Lemma \ref{secondhmsapprox} is an
algorithm which produces approximating quadruples from covering
approximations; the idea for this algorithm is from
\cite{Sapozhenko}.

\subsection{Proofs of Lemmas \ref{firsthmsapprox} and \ref{secondhmsapprox}: Approximations}
\label{subsectionapproxlemmas}

We begin with a simple observation about sums of binomial
coefficients which we will draw on repeatedly (and usually without
comment) in this section and the next. If $k=o(n)$, we have
\begin{eqnarray}
\sum_{i \leq k} {n \choose i} & \leq & (1+O(k/n)){n \choose k} \nonumber \\
& \leq & (1+O(k/n))(en/k)^k \nonumber \\
& \leq & 2^{(1+o(1))k\log(n/k)}. \label{binomial}
\end{eqnarray}

\noindent {\em Proof of Lemma \ref{firsthmsapprox}: }For each $A
\in {\cal H}$ we obtain a covering approximation for $A$ by taking
$F'(A) \subseteq G$ to be a cover of minimum size of $A$ in the
graph induced by $G \cup A$ and $P'(A) \subseteq H$ to be a cover
of minimum size of $B$ in the graph induced by $H \cup B$. Note
that $P'(A) \subseteq N(F'(A))$.

By Lemma \ref{Lconn}, $F'(A)$ is $4$-linked ($A$ is $2$-linked,
$\rho(x, F'(A)) = 1$ for each $x \in A$ and $\rho(y, A) = 1$ for
each $y \in F'(A)$). By Lemma \ref{Lcor}, $|F'(A)| \leq g(1+\ln
d)/d = O(g \log d / d)$ and $|P'(A)| \leq |H|(1+\ln d)/d = O(g
\log d / d)$ (noting that $h \leq g$).

We may therefore take ${\cal V}$ to be the set of all pairs
$(F',P') \in 2^{\cal O} \times 2^{\cal E}$ with $F'$ $4$-linked
and $P' \subseteq N(F')$, and $F', P'$ both of size at most $O(g
\log d / d)$. By Lemma \ref{Tree}, there are at most
$$
M \sum_{i \leq O(g \log d / d)} 2^{O(i \log d)} = M 2^{O(g \log^2
d / d)}
$$
possibilities for $F'$ (the factor of $M$ is for the choice of a
fixed vertex in $F'$), and, given $F'$, a further
$$
\sum_{i \leq O(g \log d / d)}{|N(F')| \choose i} =  2^{O(g \log^2 d / d)}
$$
choices for $P'$ (here we are using (\ref{binomial}) and the fact
that $|N(F')| \leq dg$). The lemma follows. \qed

\medskip

\noindent {\em Proof of Lemma \ref{secondhmsapprox}: }Fix $A
\subseteq {\cal E}$. We give an algorithm which, for input
$(F',S') \in 2^{\cal O} \times 2^{\cal E}$ satisfying $F'
\subseteq G$ and $S' \supseteq A$ produces an output $(F,S) \in
2^{\cal O} \times 2^{\cal E}$ satisfying (\ref{quad1}),
(\ref{quad2}) and (\ref{quad3}).

\medskip

Fix a linear ordering $\ll$ of $V$.

\medskip

\noindent {\bf Step $1$: }If $\{u \in A\colon  d_{G \setminus
F'}(u)\geq \sqrt{d} \}
 \neq \emptyset$, pick the smallest (with respect to $\ll$) $u$ in this
set and update $F'$ by $F' \longleftarrow F' \cup N(u)$. Repeat
this until $\{u \in A\colon  d_{G \setminus F'}(u)\geq \sqrt{d} \}
= \emptyset$. Then set $F''=F'$ and $S''=S' \setminus \{u \in
{\cal E}\colon d_{{\cal O}\setminus F''}(u) \geq \sqrt{d}\}$ and
go to Step $2$.

\medskip

\noindent {\bf Step $2$: }If $\{w \in {\cal O} \setminus G\colon
d_{S''}(w) \geq \sqrt{d} \} \neq \emptyset$, pick the smallest
(with respect to $\ll$) $w$ in this set and update $S''$ by $S''
\longleftarrow S'' \setminus N(w)$. Repeat this until $\{w \in
{\cal O} \setminus G\colon  d_{S''}(w) \geq \sqrt{d} \} =
\emptyset$. Then set $S=S''$ and $F=F''\cup \{w \in {\cal O}
\colon d_S(w) \geq \sqrt{d} \}$ and stop.

\medskip

\begin{claim} \label{algoanalysisoutput}
The output of this algorithm satisfies (\ref{quad1}),
(\ref{quad2}) and (\ref{quad3}).
\end{claim}

\noindent {\em Proof: }To see that $F \subseteq G$ and $S
\supseteq A$, first observe that $F'' \subseteq G$ (since $F'
\subseteq G$, and the vertices added to $F'$ in Step $1$ are all
in $G$) and that $S'' \supseteq A$ (or Step $1$ would not have
terminated). We then have $S \supseteq A$ since Step $2$ deletes
from $S''$ only neighbours of ${\cal O} \setminus G$, and $F
\subseteq G$ since the vertices added to $F''$ at the end of Step
$2$ are all in $G$ (or Step $2$ would not have terminated).

To verify (\ref{quad2}) and (\ref{quad3}), note that
$d_{F''}(u)>d-\sqrt{d}$ for all $u \in S''$ by definition,
$S\subseteq S''$, and $F \supseteq F''$, so that
$d_{F}(u)>d-\sqrt{d}$ for all $u \in S$; and if $w \in {\cal O}
\setminus F$ then $d_S(w) < \sqrt{d}$ (again by definition), so
that $d_{{\cal E} \setminus S}(w)
> d-\sqrt{d}$ for all $w \in {\cal O} \setminus F$. \qed

The proof of Lemma \ref{secondhmsapprox} involves a two-stage
procedure. Stage $1$ runs the algorithm described above with
$(F',{\cal E})$ as input. Stage $2$ runs it with $(P',{\cal O})$
as input and with the roles of ${\cal E}$ and ${\cal O}$ reversed.
By Claim \ref{algoanalysisoutput}, the quadruple $(F,S,P,Q)$,
where $(F,S)$ is the output of Stage $1$ and $(P,Q)$ the output of
Stage $2$, is an approximating quadruple for $A$.

\begin{claim} \label{algoanalysisnum}
The procedure described above has at most $2^{O(g \log d /
\sqrt{d})}$ outputs as the input runs over those $A \in {\cal H}$
for which $(F',P')$ is a covering approximation.
\end{claim}

\noindent Taking ${\cal W}$ to be the set of all possible outputs of the
algorithm, Lemma \ref{secondhmsapprox} follows.

\medskip

\noindent {\em Proof of Claim \ref{algoanalysisnum}: } The output
of Stage $1$ of the algorithm is determined by the set of $u$'s
whose neighbourhoods are added to $F'$ in Step $1$, and the set of
$w$'s whose neighbourhoods are removed from $S''$ in Step $2$.

Each iteration in Step $1$ removes at least $\sqrt{d}$ vertices
from $G$, so there are at most $g /\sqrt{d}$ iterations. The $u$'s
in Step $1$ are all drawn from $A$ and hence $N(F')$, a set of
size at most $d g$. So the total number of outputs for Step $1$ is
at most
$$
\sum_{i \leq g/\sqrt{d}} {d g
\choose i} = 2^{O(g \log d / \sqrt{d})}.
$$

We perform a similar analysis on Step $2$. Each $u \in
S''\setminus A$ contributes more than $d-\sqrt{d}$ edges to
$\nabla(G)$, so initially $|S''\setminus A|\leq g d/(d-\sqrt{d}) =
O(g)$. Each $w$ used in Step $2$ reduces this by at least
$\sqrt{d}$, so there are at most $O(g/\sqrt{d})$ iterations. Each
$w$ is drawn from $N(S'')$, a set which is contained in the fourth
neighbourhood of $F'$ ($S'' \subseteq N(G)$ by construction of
$S''$, $G=N(A)$ and $A \subseteq N(F')$) and so has size at most
$d^4 g$. So as with Step $1$, the total number of outputs for Step
$2$, and hence for Stage $1$, is $2^{O(g \log d / \sqrt{d})}$.

Noting that $h \leq g$, a similar analysis applied to Stage $2$
gives that that stage also has at most $2^{O(g \log d /
\sqrt{d})}$ outputs, and the claim follows. \qed

\subsection{Proof of Lemma \ref{countingas0}: Reconstruction}
\label{subsectionprovingmainapprox}

We first note an important property of approximating quadruples.
\begin{lemma} \label{sleqfqleqp}
If $(F,S,P,Q)$ is an approximating quadruple for $A \in {\cal H}$
then
\begin{eqnarray}
|S| & \leq & |F| + O(g/\sqrt{d}) \label{boundingsbyfspecific} \\
|Q| & \leq & |P| + O(h/\sqrt{d}) \label{boundingqbypspecific}.
\end{eqnarray}
\end{lemma}

\noindent {\em Proof: } Observe that $|\nabla(S,G)|$ is bounded
above by $d|F| + \sqrt{d}|G \setminus F|$ and below by $d|A| +
(d-\sqrt{d})|S \setminus A| = d|S| - \sqrt{d}|S \setminus A|$,
giving
$$
|S| \leq |F| + |(G \setminus F) \cup (S \setminus A)|/\sqrt{d},
$$
and that each $u \in (G \setminus F) \cup (S \setminus A)$
contributes more than $d-\sqrt{d}$ edges to $\nabla(G)$, a set of
size $g d$, giving
$$
|(G \setminus F) \cup (S \setminus A)| \leq 2 g d/(d-\sqrt{d}) = O(g).
$$
These two observations together give (\ref{boundingsbyfspecific}).
The proof of (\ref{boundingqbypspecific}) is similar. \qed

Lemma \ref{countingas0} now follows from
\begin{lemma} \label{countingas}
For each $(F,S,P,Q) \in 2^{{\cal O}} \times 2^{{\cal E}} \times
2^{{\cal E}} \times 2^{{\cal O}}$ satisfying
(\ref{boundingsbyfspecific}) and (\ref{boundingqbypspecific}),
there are at most $2^{g - b - \Omega(g/\log d)}$ $A$'s in ${\cal
H}$ satisfying (\ref{quad1}) and (\ref{quad4}).
\end{lemma}

\noindent {\em Proof: } For $A \in {\cal H}$, write
$$
[A]=\{u \in {\cal E}\colon N(u)\subseteq N(A)\},
$$
and write $a'$ for $|[A]|$. Note
that although $G$ does not determine $A$, it does determine $[A]$.
By Lemma \ref{boundsondelta}, there is an absolute constant
$\gamma>0$ (independent of $a', g, b$ and $h$) such that
\beq{finallyusingdelta} g-a' > \gamma g ~~~~~ \mbox{and} ~~~~~ h-b
> \gamma h. \enq

Say that $Q$ is {\bf tight} if $|Q| < b + \gamma h/\log d$, and
{\bf slack} otherwise, and that $S$ is {\bf tight} if $|S| < g -
\gamma g/(4\log d)$ and {\bf slack} otherwise.

We now describe a procedure which, for input $(F,S,P,Q)$, produces
an output $A$ which satisfies (\ref{quad1}) and (\ref{quad4}). The
procedure involves a sequence of choices, the nature of the
choices depending on whether $S$ and $Q$ are tight or slack.

We begin by identifying a subset $D$ of $A$ which can be specified
relatively ``cheaply'': if $Q$ is tight, we pick $B \subseteq Q$
with $|B|=b$ and take $D=N(B)$; if $Q$ is slack, we simply take
$D=P$ (recalling that $P \subseteq H \subseteq A$).

If $S$ is tight, we complete the specification of $A$ by choosing
$A \setminus D \subseteq S \setminus D$. If $S$ is slack, we first
complete the specification of $G$ by choosing $G \setminus F
\subseteq N(S) \setminus F$. Note that in this case,
(\ref{boundingsbyfspecific}) implies \beq{glessf} |G \setminus F|
< \gamma g /(3\log d). \enq We then complete the specification of
$A$ by choosing $A \setminus D \subseteq [A] \setminus D$ (noting
that we do know $[A] \setminus D$ at this point).

This procedure produces all possible $A \in {\cal H}$ satisfying
(\ref{quad1}) and (\ref{quad4}) (and more). Before bounding the
number of outputs, we gather together some useful observations.

From (\ref{boundingsbyfspecific}) and (\ref{boundingqbypspecific}) we have
\beq{sq}
|S| = O(g) ~~~~~ \mbox{and} ~~~~~ |Q| = O(h).
\enq

If $Q$ is tight then there are at most
\begin{eqnarray}
\sum_{i \leq \gamma h /\log d}{|Q| \choose |Q|-i} & \leq &
           \sum_{i \leq \gamma h /\log d}{O(h) \choose i} \nonumber \\
                                                    & \leq &
           2^{O(\gamma h /\log d)\log(O(\log d/\gamma))} \nonumber \\
                                                    & \leq &
           2^{\gamma h/2} \label{choicesfordqsmall}
\end{eqnarray}
possibilities for $D$, and in this case $|D|=h$; while if $Q$ is
slack there is just one possibility for $D$, and in this case
(using (\ref{boundingqbypspecific}))
\begin{eqnarray}
|D|=|P| & > & |Q|-\Omega(h/ \sqrt{d}) \nonumber \\
        & > & b+\gamma h /\log d  -\Omega(h/\sqrt{d}) \nonumber \\
        & \geq  & b+\gamma h /(2\log d). \label{qlarge}
\end{eqnarray}

If $S$ is slack then (since $|N(S)\setminus F| \leq d|S| \leq
O(dg)$; see (\ref{sq})) the number of possibilities for $G
\setminus F$ is at most
\begin{eqnarray}
\sum_{i<\gamma g /(3\log d)}{O(gd) \choose i} & \leq &
2^{(1+o(1))(\gamma g /(3\log d))\log(O(d\log d/\gamma))} \nonumber \\
 & \leq & 2^{\gamma g/2}. \label{choicesforglessfslarge}
\end{eqnarray}

We now bound the number of outputs of the procedure, considering
separately the four cases determined by whether $S$ and $Q$ are
slack or tight.

If $S$ and $Q$ are both tight then the number of possibilities for
$A$ is at most \beq{ssmallqsmall} 2^{[\gamma h/2]+[g - \gamma
g/(4\log d)-h]}
      < 2^{g - \gamma g/(4\log d)-b-\gamma h/2}.
\enq (The first term in the exponent on the left-hand side
corresponds to the choice of $D$ (using
(\ref{choicesfordqsmall})), and the second to the choice of
$A\setminus D$ (note that since $S$ and $Q$ are both tight,
$|S\setminus D| \leq g - \gamma g/(4\log d)-h$). To get the
right-hand side, we use the second part of
(\ref{finallyusingdelta}).)

If $S$ is tight and $Q$ is slack then the total is at most
\beq{ssmallqlarge} 2^{[g - \gamma g/(4\log d)-b-\gamma h/(2\log
d)]}. \enq (Here there is no choice for $D$, and the exponent
corresponds to the choice of $A\setminus D$ (using
(\ref{qlarge})).)

If $Q$ is tight then $|[A]\setminus D| = a'-h$, so that if $S$ is
slack (and $Q$ tight) then the number of possibilities for $A$ is
at most
\begin{equation} \label{slargeqsmall}
2^{[\gamma h/2]+[\gamma g/2]+[a'-h]}<2^{g-\gamma g/2-b-\gamma
h/2}.
\end{equation}
(The first term on the left-hand side corresponds to the choice of
$D$ (using (\ref{choicesfordqsmall}), the second to the choice of
$G\setminus F$ (using (\ref{choicesforglessfslarge})) and the
third to the choice of $A\setminus D$. On the right-hand side, we
use both parts of (\ref{finallyusingdelta}).)

Finally, if $Q$ is slack then $|[A]\setminus D| \leq a'- b-\gamma
h /(2\log d)$ (see (\ref{qlarge})), so that if $S$ and $Q$ are
both slack the number of possibilities for $A$ is at most
\begin{equation} \label{slargeqlarge}
2^{[\gamma g/2]+[a'-b-\gamma h/(2\log d)]}
     < 2^{g-\gamma g/2-b-\gamma h/(2\log d)}.
\end{equation}
(The first term on the left-hand side corresponds to the choice of
$G\setminus F$ and the second to the choice of $A\setminus D$. The
right-hand side uses the first part of (\ref{finallyusingdelta}).)

Noting that $h \leq g$, the lemma follows from (\ref{ssmallqsmall}), (\ref{ssmallqlarge}),
(\ref{slargeqsmall}) and (\ref{slargeqlarge}). \qed

\section{Proof of Theorem \ref{heartofmatter}} \label{sectionprovingo(1)}

We say that a nice $A \subseteq {\cal E}$ is {\bf of type I} if
$|A|<d/2$, {\bf of type II} if $d/2 \leq |A|< d^{2}$ and {\bf of
type III} otherwise. We consider the portions of the sum in
Theorem \ref{heartofmatter} corresponding to type I, II and III
$A$'s separately.

If $A$ is of type I, then by Lemma \ref{boundsondeltaAsmall},
$|N(A)| \geq d|A| - 2|A|(|A|-1)$. Note also that in this case,
$B(A)=\emptyset$. By Lemma \ref{Tree}, for each $2\leq i < d/2$,
there are at most $M2^{O(i\log d)} < 2^{d+O(i\log d)}$ $2$-linked
subsets of ${\cal E}$ of size $i$. So
\begin{eqnarray}
\sum_{\mbox{$A$ of type I}} 2^{-|N(A)|+|B(A)|} & \leq &
\sum_{i=2}^{d/2} 2^{d+O(i\log d)-di+2i(i-1)} \nonumber \\
& = & e^{-\Omega(d)}. \label{typeIA}
\end{eqnarray}

We do something similar if $A$ is of type II. Here Lemma
\ref{boundsondeltaAsmall} gives $|N(A)| \geq \Omega(d)|A|$ and $|B(A)| \leq O(1/d)|A|$
(recalling that $N(B) \subseteq A$), and so
\begin{eqnarray}
\sum_{\mbox{$A$ of type II}} 2^{-|N(A)|+|B(A)|} & \leq &
\sum_{i=d/2}^{d^{2}} 2^{d+O(i\log d)-\Omega(d)i+O(1/d)i} \nonumber \\
& = & e^{-\Omega(d)}. \label{typeIIA}
\end{eqnarray}

We partition the set of $A$'s of type III according to the sizes
of $A$, $N(A), B(A)$ and $H(A)~(=N(B(A)))$ and use Lemma
\ref{mainapprox} to bound the sizes of the partition classes. In
this case we have $|N(A)| \geq d^{2}$. So (summing only over those
values of $a$, $g$, $b$ and $h$ for which ${\cal H}(a,g,b,h) \neq
\emptyset$ and $g \geq d^2$, and with the inequalities justified
below)
\begin{eqnarray}
\sum_{\mbox{$A$ of type III}} 2^{-|N(A)|+|B(A)|} & = & \sum_{a, g, b, h} |{\cal H}(a,g,b,h)| 2^{-g+b} \nonumber \\
  & \leq & M \sum_{a, g, b, h} 2^{-\Omega(g/\log d)} \label{usingmainapprox} \\
  &   <  & M^4 \sum_{g \geq d^{2}} 2^{-\Omega(g/ \log d)} \label{msquared} \\
  & \leq & \left(M^4/(1-2^{-\Omega(1/\log d)})^2\right) 2^{-\Omega(d^{2}/\log d)} \nonumber \\
  &   =  & e^{-\Omega(d)}. \label{typeIIIAgpart}
\end{eqnarray}
Here (\ref{usingmainapprox}) is from Lemma \ref{mainapprox} and in (\ref{msquared}) we use the fact that there
are fewer than $M$ choices for each of $a$, $b$ and $h$.

Combining (\ref{typeIA}), (\ref{typeIIA}) and
(\ref{typeIIIAgpart}), we have Theorem \ref{heartofmatter}. \qed

\section{Remarks} \label{sectionremarks}

The point of departure for our proof of Theorem \ref{ourresult} is
Lemma \ref{entropylemma}, which allows us to focus immediately on
those homomorphisms which are predominantly single-valued on one
side of the cube. The proof of this lemma given in \cite{Kahn}
relies heavily on the structure of the cube (in particular on the
fact that the neighbourhoods of adjacent vertices induce a perfect
matching), and it does not seem obvious at the moment how to get
beyond this and generalize Theorem \ref{ourresult} to a larger
class of graphs.

On the other hand, the proofs of Theorem \ref{heartofmatter} and
Lemma \ref{mainapprox} are much less dependent on the specific
structure of the cube, using only the isoperimetric bounds of
Section \ref{sectionisoperimetry}. As such, it should be possible
to extend these results considerably. To illustrate this, it is
worth comparing Lemma \ref{mainapprox} with the main lemma of
\cite{Sapozhenko}. To state that, we need some notation. Let $G$
be a $d$-regular bipartite graph with bounded co-degree (i.e.,
every pair of vertices has a bounded number of common neighbours).
Write $X$ and $Y$ for the bipartition classes of $G$. For any $a'$
and $g$, set
$$
{\cal G}(a',g)=\{A \subseteq X \colon \mbox{$A$ $2$-linked,
$|N(A)|=g, |[A]|\leq a'$}\},
$$
(recall that $[A]=\{x\in X \colon N(x)\subseteq N(A)\}$), and set
$\delta=(g-a')/g$. Using slightly more versatile notions of
approximation than those introduced in Section
\ref{sectionmainapproximation}, the following is proved in
\cite{Sapozhenko}:
\begin{thm} \label{Sapslemma}
For $d$ sufficiently large, and for any $a'$ and $g$ satisfying
$1> \delta >\log^9 d/d^2$,
$$
|{\cal G}(a',g)| \leq |X|2^{g(1-\delta/(6\log d))}.
$$
\end{thm}

Notice that (by the results of Section \ref{sectionisoperimetry})
the sum in Theorem \ref{heartofmatter} is extending only over sets
$A$ which satisfy $(|N(A)|-|A|)/|N(A)| \geq \Omega(1)$, a much
stronger condition than that imposed in Theorem \ref{Sapslemma}.
By slightly modifying our notions of approximation, we may extend
the validity of Lemma \ref{mainapprox} to cover a similar range as
Theorem \ref{Sapslemma}. However, the analysis is considerably
more involved, and we do not do so here.

\bigskip

\noindent {\bf Acknowledgement: }The author thanks Jeff Kahn for numerous helpful discussions.


\begin{thebibliography}{99}

\bibitem{Athanasiadis}
C. A. Athanasiadis,
{\em Algebraic combinatorics of graph spectra, subspace arrangements,
and Tutte polynomials}, thesis, Massachusetts Institute of Technology, 1996.

\bibitem{BenjaminiHaggstromMossel}
I. Benjamini, O. H\"aggstr\"om and E. Mossel, On random graph
homomorphisms into ${\bf Z}$, {\em J. Combinatorial Th. (B)} {\bf
78} no. 1 (2000), 86--114.

\bibitem{Bollobas4}
B. Bollob\'as,
{\em Combinatorics}, Cambridge University Press, Cambridge, 1986.

\bibitem{Bollobas}
B. Bollob\'as,
{\em Modern Graph Theory}, Springer, New York, 1998.

\bibitem{Diestel}
R. Diestel, {\em Graph Theory}, Springer, New York, 1997.

\bibitem{Furedi}
Z. F\"uredi, Matchings and covers in hypergraphs,
{\em Graphs and Comb.} {\bf 4} (1988), 115--206.

\bibitem{Kahn}
J. Kahn, Range of the cube-indexed random walk, {\em Israel J.
Math.} {\bf 124} (2001) 189--201.

\bibitem{KahnLawrenz}
J. Kahn and A. Lawrenz,
Generalized rank functions and an entropy argument,
{\em J. Combinatorial Th. (A)} {\bf 87} (1999), 398--403.

\bibitem{Knuth}
D. Knuth, {\em The Art of Computer Programming} Vol. I,
Addison Wesley, London, 1969.

\bibitem{KornerWei}
J. K\"orner and V. Wei, Odd and even Hamming spheres also have
minimum boundary, {\em Discrete Math.} {\bf 51} (1984), 147--165.

\bibitem{KorshunovSapozhenko}
A. D. Korshunov and A. A. Sapozhenko,
The number of binary codes with distance $2$,
{\em Problemy Kibernet.} {\bf 40} (1983), 111--130. (Russian)

\bibitem{Lovasz}
L. Lov\'asz, On the ratio of optimal integral and fractional
covers, {\em Discrete Math.} {\bf 13} (1975) 383--390.

\bibitem{Randall2}
D. Randall, personal communication.

\bibitem{Sapozhenko}
A. A. Sapozhenko,
On the number of connected subsets with given cardinality of the boundary in bipartite graphs,
{\em Metody Diskret. Analiz.} {\bf 45} (1987), 42--70.  (Russian)

\bibitem{Sapozhenko2}
A. A. Sapozhenko, The number of antichains in ranked partially
ordered sets, {\em Diskret. Mat.} {\bf 1} (1989), 74--93.
(Russian; translation in Discrete Math. Appl. 1 no. 1 (1991),
35--58)

\bibitem{Stein}
S. K. Stein, Two combinatorial covering theorems, {\em J.
Combinatorial Th. (A)} {\bf 16} (1974), 391--397.

\end{thebibliography}
\end{document}